\documentclass[12pt]{article}
\usepackage{latexsym,amssymb,amsfonts}
\usepackage{mathrsfs}
\usepackage{psfrag}
\usepackage{amsmath, amsbsy}
\usepackage{amsopn, amstext}

\def\sqr#1#2{{\vcenter{\vbox{\hrule height.#2pt
              \hbox{\vrule width.#2pt height#1pt \kern#1pt \vrule width.#2pt}
              \hrule height.#2pt}}}}
\def\signed #1{{\unskip\nobreak\hfil\penalty50
              \hskip2em\hbox{}\nobreak\hfil#1
              \parfillskip=0pt \finalhyphendemerits=0 \par}}
\def\endpf{\signed {$\sqr69$}}

\def\dbR{{\mathop{\rm l\negthinspace R}}}

\def\1n{\negthinspace }
\def\2n{\negthinspace \negthinspace }
\def\3n{\negthinspace \negthinspace \negthinspace }

\def\ns{\noalign{\ss}}

\def\dbF{{\mathbb{F}}}

\def\dbH{{\mathbb{H}}}

\def\dbP{{\mathbb{P}}}

\def\dbR{{\mathbb{R}}}

\def\pa{\partial}

\def\cd{\cdot}
\def\cds{\cdots}

\def\as{\hbox{\rm a.s.{ }}}

\def\({\Big (}
\def\){\Big )}
\def\[{\Big[}
\def\]{\Big]}

\def\={\buildrel \triangle \over =}

\def\-{\mbox{-}}
\def\={\buildrel \triangle \over =}

\def\ds{\displaystyle}
\def\mE{\mathbb{E}}
\def\a{\alpha}

\def\g{\gamma}
\def\d{\delta}
\def\e{\varepsilon}

 \def\n{\nabla}
\def\si{\sigma}
\def\t{\times}

\def\th{\theta}
\def\om{\omega}

\def\D{\Delta}

\def\Om{\Omega}

\def\cA{{\cal A}}

\def\cD{{\cal D}}

\def\cM{{\cal M}}

\def\cP{{\cal P}}

\def\cV{{\cal V}}

\def\no{\noindent}

\def\ss{\smallskip}
\def\ms{\medskip}
\def\bs{\bigskip}
\def\q{\quad}
\def\qq{\qquad}

\def\pa{\partial}

\def\cd{\cdot}
\def\cds{\cdots}

\def\as{\hbox{\rm a.s.}}

\def\|{\Big |}
\def\({\Big (}
\def\){\Big )}
\def\[{\Big[}
\def\]{\Big]}

\def\Om{\Omega}

\def\max{\mathop{\rm max}}
\def\min{\mathop{\rm min}}

\def\pa{\partial}

\def\cd{\cdot}
\def\cds{\cdots}

\def\as{\hbox{\rm a.s.{ }}}

\def\|{\Big |}
\def\({\Big (}
\def\){\Big )}
\def\[{\Big[}
\def\]{\Big]}

\newtheorem{lemma}{Lemma}[section]
\newtheorem{remark}{Remark}[section]

\newtheorem{theorem}{Theorem}[section]
\newtheorem{corollary}{Corollary}[section]

\newtheorem{definition}{Definition}[section]

\makeatletter
   
   \@addtoreset{equation}{section}
\makeatother

\allowdisplaybreaks

\begin{document}
\title{\bf Unique Continuation for Stochastic Hyperbolic Equations}

\author{Qi L\"u\thanks{School of Mathematics,
Sichuan University, Chengdu 610064, Sichuan
Province, China. The author is partially
supported by NSF of China under grants 11471231,
the Fundamental Research Funds for the Central
Universities in China under grant 2015SCU04A02
and Grant MTM2011-29306-C02-00 of the MICINN,
Spain. {\small\it E-mail:} {\small\tt
lu@scu.edu.cn}.} \q and\q Zhongqi Yin\thanks{Department of Mathematics, Sichuan Normal University, Chengdu 610068, P. R. China. The author is supported by the General Fund Project of Sichuan Provincial Department of Education in China under grant 13ZB0164. {\small\it E-mail:} {\small\tt zhongqiyin@sicnu.edu.cn}.}}

\date{}

\maketitle

\begin{abstract}\no
In this paper, we derive a local unique
continuation property for stochastic hyperbolic
equations without boundary conditions. This
result is proved by a global Carleman estimate.
As far as we know, this is the first result in
this topic.
\end{abstract}

\bs

\no{\bf 2000 Mathematics Subject
Classification}.  Primary   60H15, 93B07.

\bs

\no{\bf Key Words}. Stochastic hyperbolic equation, Carleman
estimate, unique continuation.

\ms

\section{Introduction }
\q Let $T > 0$, $G \subset \mathbb{R}^{n}$ ($n \in
\mathbb{N}$) be a given bounded domain.
Throughout this paper, we will use $C$ to denote
a generic positive constant depending  only on
$T$ and $G$, which may change from line to line.
Denote $Q = G\times (0, T)$.

Let $(\Om, {\cal F}, \dbF, \dbP)$ with
$\dbF\=\{{\cal F}_t\}_{t \geq 0}$ be a complete
filtered probability space on which a  one
dimensional standard Brownian motion
$\{W(t)\}_{t\geq 0}$ is defined. Assume that $H$
is a Fr\'echet space. Let $L^2_{\dbF}(0, T; H)$
be the Fr\'echet space consisting all $H$-valued
$\dbF$-adapted process $X(\cdot)$ such that $\mE
|X(\cdot)|_{L^2_{\dbF}(0, T; H)}^2 < +\infty$,
on which the canonical quasi-norm is endowed. By
$L^{\infty}_{\dbF}(0, T; H)$ we denote the
Fr\'{e}chet space of all $H$-valued
$\dbF$-adapted bounded processes equipped with
the canonical quasi-norm and by
$L^2_{\dbF}(\Omega; C([0, T]; H))$ the
Fr\'{e}chet space of all $H$-valued
$\dbF$-adapted continuous precesses $X(\cdot)$
with $\mE |X(\cdot)|^2_{C([0, T]; H)} < +\infty$
and equipped with the canonical quasi-norm.

This paper is devoted to the study of a local
unique continuation property for the following
stochastic hyperbolic equation:
\begin{equation}\label{system1}
\sigma dz_{t} - \D z dt = \big(b_1 z_t +
b_2\cd\nabla z + b_3 z \big)dt +  b_4 z dW(t)
\qq {\mbox { in }} Q,
\end{equation}
where, $\sigma\in C^1(\overline Q)$ is positive,
\begin{eqnarray*} \label{aibi}
&\,& b_1 \in
L_{\dbF}^{\infty}(0,T;L^{\infty}_{loc}(G)), \qq
b_2
\in L_{\dbF}^{\infty}(0,T;L^{\infty}_{loc}(G;\mathbb{R}^{n})), \nonumber \\
&\,& b_3 \in
L_{\dbF}^{\infty}(0,T;L^{n}_{loc}(G)), \qq \,\,
b_4 \in
L_{\dbF}^{\infty}(0,T;L^{\infty}_{loc}(G)).
\end{eqnarray*}
Put
\begin{equation}\label{HT}
\dbH_{T} \= L_{\dbF}^2 (\Om;
C([0,T];H_{0,loc}^1(G)))\cap L_{\dbF}^2 (\Om;
C^{1}([0,T];L^2_{loc}(G))).
\end{equation}

The definition of  solutions to equation
(\ref{system1}) is given in the following sense.

\begin{definition} \label{def solution to sys} We call $z
\in \dbH_T$  a solution of equation
(\ref{system1}) if for each $t \in [0,T]$,
$G'\subset\subset G$ and  $\eta \in H_0^1(G')$,
it holds that
\begin{equation}\label{solution to sysh}
\begin{array}{ll}\ds
\q \int_{G'} z_t(t,x)\eta(x)dx - \int_{G'}
z_t(0,x)\eta(x)dx -\int_0^t \int_{G'} \si_t(t,x)z_t(t,x)\eta(x)dx\\
\ns\ds = \int_0^t \int_{G'}  \[ -\nabla
z(s,x)\cd\nabla\eta(x) + \big(b_1 z_t +
b_2\cd\nabla z + b_3 z\big)\eta(x) \]dxds \\
\ns\ds \q + \int_0^t \int_{G'}   b_4 z \eta(x)
dx dW(s), \qq \dbP\mbox{-a.s. }
\end{array}
\end{equation}
\end{definition}

Let $S \subset\subset  G$ be a
$C^2$-hypersurface. Let $x_0\in S\setminus
\partial G$ and suppose that $S$
divides the ball $B_{\rho}(x_0)\subset G$, centered at $x_0$ and with radius $\rho$, into
two parts $\cD_{\rho}^+$ and $\cD_{\rho}^-$. Denote
as usual by $\nu(x)$ the unit normal vector to
$S$ at $x$ inward to $\cD_{\rho}^+$.

Let $y$ be a solution to equation
(\ref{system1}). Let $\e>0$. This paper is
devoted to the following local unique
continuation problem:

\medskip

{\bf (Pu)} Can $y$ in $\cD_{\rho}^+\times
(\e,T-\e)$ be uniquely determined by the values
of $y$ in $\cD_{\rho}^-\times (0,T)$?

\medskip

In other words, Problem {\bf (Pu)} concerns that
whether the values of the solution in one side of
$S$ uniquely determine its values in the other
side. Clearly, it is equivalent to the following
problem:

\medskip

{\bf (Pu1)} Can we conclude that $y=0$ in
$\cD_{\rho}^+\times (\e,T-\e)$, provided that
$y=0$ in $\cD_{\rho}^-\times (0,T)$?

\medskip

Unique continuation problems for deterministic
PDEs are studied extensively in the literature.
Generally speaking, a unique continuation result
is a statement in the following sense:

\vspace{0.2cm}

Let $u$ be a solution of a PDE and two regions
$\cM_1\subset\cM_2$. Then $u$ is determined uniquely
in $\cM_2$ by its values in $\cM_1$.

\vspace{0.2cm}

Problem {\bf (Pu)} is a natural generalization
of the unique continuation problems under the
stochastic setting, i.e., $\cM_1=\cD^-_\rho$ and
$\cM_2=B_\rho (x_0)$.

There is a long history for the study of unique
continuation property (UCP for short) for
deterministic PDEs. Classical results date back
to Cauchy-Kovalevskaya theorem and Holmgren's
uniqueness theorem. These two results need the
coefficients of the PDE  analytic to get
the UCP. In 1939,  T. Carleman introduced in
\cite{Carleman1} a new method, which was based
on weighted estimates, to prove UCP for two
dimensional elliptic equations with $L^\infty$
coefficients.  This method, which is called
``Carleman estimates method" nowadays,  turned
out to be quite powerful and has been
developed extensively in the literature and
becomes the most useful tool to obtain UCP for
PDEs (e.g.
\cite{Calderon,Hormander2,Kenig,Tataru,Tataru1,Zuily}).
In particular,  unique continuation results for
solutions of hyperbolic equations across
hypersurfaces were studied by many authors (e.g.
\cite{Hormander3,Kenig1,Robbiano1,Sogge,Tataru2}).

Compared with the deterministic PDEs, there are
very few results concerning UCP for stochastic
PDEs. To our best knowledge, \cite{Zh} is the
first result for UCP of stochastic PDEs, in
which the  author shows that  a solution to a
stochastic parabolic equation vanishes, provided
that it vanishes in any subdomain and this
result was improved in \cite{LL,LY} where less
geometric condition is assumed to the set where
the solution vanishes. The first result of UCP
for stochastic hyperbolic equations was obtained
in \cite{Zhangxu3}. Some improvements were made in
\cite{Lu,LZ2}. The results in  \cite{Lu}
and \cite{Zhangxu3} concerned the global UCP for
stochastic hyperbolic equations with a
homogeneous Dirichlet boundary condition, i.e.,
they concluded that the solution to a stochastic
hyperbolic equation vanishes, provided that it
equals zero in a large enough subdomain. In this
paper, we focus on the local UCP for stochastic
hyperbolic equations without boundary condition,
that is, can a solution be determined locally?

To present the main result of this paper, let us
first introduce the following notion.

\begin{definition}\label{def1}
Let $x_0\in S$ and $ K>0$. $S$ is said to
satisfy the outer paraboloid condition with $K$
at $x_0$ if there exists a neighborhood $\cV$ of
$x_0$ and a paraboloid $\cP$ tangential to $S$
at $x_0$ and $\cP\cap\cV\subset  \cD^{-}_\rho$
with $\cP$ congruent to $\ds x_1 = K\sum_{j=2}^n
x_j^2$.
\end{definition}

The main result in this paper is the following
one.

\begin{theorem}\label{th1}
Let $x_0\in S\setminus\partial S$ such that
$\frac{\pa \sigma(x_0,T/2)}{\pa\nu}<0$, and
let $S$ satisfy the outer paraboloid
condition with
\begin{equation}\label{th1-eq1}
K <\frac{-\frac{ \partial
\sigma}{\pa\nu}(x_0,T/2)}{4(|\sigma|_{L^\infty(B_\rho(x_0,T/2))}+1)}.
\end{equation}
Let $z\in \dbH_T$ be a solution of the equation \eqref{system1}
satisfying that
\begin{equation}\label{th1-eq2}
z=\frac{\pa z}{\partial\nu} = 0 \q\mbox{ on }
(0,T)\times S,\;\;\dbP\mbox{-}\as
\end{equation}
Then, there is a neighborhood $\cV$ of $x_0$ and
$\e\in (0,T/2)$ such that
\begin{equation}\label{th1-eq3}
z = 0 \q\mbox{ in } (\cV\cap \cD_{\rho}^+)\times
(\e,T-\e),\;\;\dbP\mbox{-}\as
\end{equation}
\end{theorem}
\begin{remark}
In Theorem \ref{th1}, we assume that $\frac{\pa
\sigma(x_0,T/2)}{\pa\nu}<0$. It is related to the
propagation of wave. This is a reasonable
assumption since the UCP may not hold if it is
not fulfilled (e.g. \cite{Alinhac}). It can be
regarded as a kind of pseudoconvex condition
(e.g. \cite[Chapter XXVII]{Hormander2}).
\end{remark}
\begin{remark}
If $S$ is a hyperplane, then  Condition
\ref{th1-eq1} always satisfies since we can take
$K = 0$.
\end{remark}
\begin{remark}
From Theorem \ref{th1}, one can get many
classical UCP results for deterministic
hyperbolic equations (e.g.
\cite{Hormander3,Robbiano,Tataru2}).  
\end{remark}

As an immediate corollary of Theorem \ref{th1},
we have the following UCP.

\begin{corollary}\label{cor1}
Let $x_0\in S \setminus\partial S$ such that
$\frac{\partial \sigma(x_0,T/2)}{\pa\nu}<0$, and let
$S$ satisfy the outer paraboloid condition with
\begin{equation}\label{th1-eq1-1}
K <\frac{-\frac{ \partial
\sigma}{\pa\nu}(x_0,T/2)}{\ds 4(|\sigma|_{L^\infty(B_\rho(x_0,T/2))}+1)}.
\end{equation}
Then   for any $z\in \dbH_T$ which solves
equation \eqref{system1} satisfying that
\begin{equation}\label{th1-eq2}
z= 0 \q\mbox{ on } \cD^-_\rho\times
(0,T),\;\dbP\mbox{-}\as,
\end{equation}
there is a neighborhood $\cV$ of $x_0$ and
$\e\in (0,T/2)$ such that
\begin{equation}\label{th1-eq3}
z = 0 \q\mbox{ in }
\big(\cV\cap\cD^+_\rho\big)\times
(\e,T-\e),\;\dbP\mbox{-}\as
\end{equation}
\end{corollary}

Similar to the deterministic settings, we shall
use the stochastic versions of Carleman estimate
for stochastic hyperbolic equations to establish
our estimate. Carleman estimates for stochastic
PDEs are studied extensively in recent years
(see
\cite{barbu1,Liuxu1,Luqi3,Luqi7,Lu,Tang-Zhang1,Zhangxu6}
and the reference therein). Carleman estimate
for stochastic hyperbolic equations was first
obtained in \cite{Zhangxu3}. Compared with the
result in \cite{Zhangxu3}, we need to handle a
more complex case (see Section \ref{sec-point}
for more details).

The rest of this paper is organized as follows.
In Section \ref{sec-point}, we derive a
point-wise estimate for stochastic hyperbolic
operator, which is crucial for us to establish the
desired Carleman estimate in this paper. In
Section \ref{sec-weight}, we explain the choice
of weight function in the Carleman estimate. Section \ref{sec-car} is devoted to the proof of
a Carleman estimate while Section \ref{sec-main}
is addressed to the proof of the main result.


\section{A point-wise estimate for stochastic hyperbolic
operator}\label{sec-point}


We introduce the
following point-wise Carleman estimate for stochastic hyperbolic operators. This estimate has its own independent interest and  will be of particular importance in the proof for the main result.
\begin{lemma}\label{hyperbolic1}
Let $\ell,\Psi \in C^2((0,T)\t\mathbb{R}^n)$.
Assume $u$ is an
$H^2_{loc}(\mathbb{R}^n)$-valued $\dbF$-adapted
process such that $u_t$ is an
$L^2(\mathbb{R}^n)$-valued semimartingale. Set
$\theta = e^\ell$ and $v=\theta u$. Then, for
a.e. $x\in \mathbb{R}^n$ and $\dbP$-a.s. $\om
\in \Om$, %
\begin{eqnarray}\label{hyperbolic2} \ds
& \ds  \theta \( -2\sigma\ell_t v_t +
2\sum_{i,j=1}^n b^{ij}\ell_i v_j + \Psi v \)
\[ \sigma du_t - \sum_{i,j=1}^n (b^{ij}u_i)_j dt \]\nonumber \\
 & \ds+\sum_{i,j=1}^n \[ \sum_{i',j'=1}^n \(
2b^{ij}b^{i'j'}\ell_{i'}v_iv_{j'} -
b^{ij}b^{i'j'}\ell_i v_{i'}v_{j'}
\) - 2b^{ij}\ell_t v_i v_t + \sigma b^{ij}\ell_i v_t^2 \nonumber \\
 & \ds + \Psi b^{ij}v_i v - \Big( A\ell_i +
\frac{1}{2}\Psi_i\Big)b^{ij}v^2 \]_j dt +d\[
\sigma\sum_{i,j=1}^n b^{ij}\ell_t v_i v_j \nonumber\\
  & \ds - 2\sigma\sum_{i,j=1}^n
b^{ij}\ell_iv_jv_t  + \sigma^2\ell_t v_t^2 -
\sigma\Psi v_t v + \Big( \sigma A\ell_t +
\frac{1}{2}(\sigma\Psi)_t\Big)v^2 \] \\
 \ds = &\ds \bigg\{ \[(\sigma^2\ell_{t})_t +
\sum_{i,j =1}^n (\sigma b^{ij}\ell_i)_{j} -
\sigma\Psi \]v_t^2 - 2\sum_{i,j=1}^n [(\sigma
b^{ij}\ell_j)_t +
b^{ij}(\sigma\ell_{t})_j]v_iv_t \nonumber \\
\ns& \ds +\sum_{i,j=}^n \Big[ (\sigma
b^{ij}\ell_t)_t + \sum_{i',j'=1}^n
\Big(2b^{ij'}(b^{i'j}\ell_{i'})_{j'}-(b^{ij}b^{i'j'}\ell_{i'})_{j'}\Big)
+ \Psi b^{ij} \]v_iv_j \nonumber\\
 & \ds + Bv^2 + \Big( -2\sigma\ell_tv_t +
2\sum_{i,j=1}^n b^{ij}\ell_iv_j + \Psi v
\Big)^2\bigg\} dt + \sigma^2\theta^2
\ell_t(du_t)^2,\nonumber
\end{eqnarray}
where $(du_t)^2$ denotes the quadratic
variation process of $u_t$, and $A$ and
$B$ are stated as follows:
\begin{equation}\label{AB1}
\left\{
 \begin{array}{ll}
 \ds A\=\sigma (\ell_t^2-\ell_{tt})-\sum_{i,j=1}^n
 (b^{ij}\ell_i\ell_j-b^{ij}_j\ell_i
 -b^{ij}\ell_{ij})-\Psi,\\
\ns
 \ds
 B\=A\Psi+(\sigma A\ell_t)_t-
 \sum_{i,j}(Ab^{ij}\ell_i)_j +\dfrac 12 \[(\sigma\Psi)_{tt}-\sum_{i,j=1}^n (b^{ij}\Psi_i)_j\].
 \end{array}
\right.
\end{equation}
\end{lemma}

\begin{remark}
When $\si=1$, equality \eqref{hyperbolic2} had
been established in \cite{Zhangxu3}. The
computation for the general $\si$ is more
complex. One needs to handle the terms
concerning $\si$ carefully.
\end{remark}

{\it Proof of Lemma \ref{hyperbolic1}.} By
$v(t,x)=\th(t,x)u(t,x)$, we have
$$u_t=\th^{-1}(v_t-\ell _tv), u_j=\th^{-1}(v_j-\ell _jv)$$ for
$j=1,2,\ldots,n$. Then, for that $\theta$ is deterministic, we have
\begin{equation}
\label{4.15-eq1}
\begin{array}
{ll}
   \sigma du_t =  \ds \sigma d[\theta^{-1}(v_t - \ell_t v)]= \ds \sigma\theta^{-1}\[dv_t -2\ell_t v_t dt  + (\ell_t ^2 - \ell_{tt} ) v dt\].
   \end{array}
\end{equation}
Moreover, we find that
\begin{equation}
\label{4.15-eq2}
\begin{array}
{ll}
 \ds\sum_{i,j=1}^n (b^{ij}u_i)_j \3n & = \ds  \sum_{i,j=1}^n \(b^{ij}\theta^{-1}(v_i- \ell_i v)\)_j \\
  \ns & = \ds \theta^{-1}\sum_{i,j=1}^n\[(b^{ij}v_i)_j - 2b^{ij}\ell_i v_j+( b^{ij}\ell_i \ell_j  - b^{ij}_j\ell_i  - b^{ij}\ell_{ij}) v \].
   \end{array}
\end{equation}
As an immediate result of \eqref{4.15-eq1} and \eqref{4.15-eq2}, we have that
\begin{equation}
\label{4.15-eq1-eq2}
\begin{array}
{ll}
   &  \ds\sigma d u_t - \sum_{i,j=1}^n (b^{ij}u_i)_j dt\\
     \ns  = & \3n\ds \theta^{-1}\[\(\sigma dv_t - \sum_{i,j=1}^n (b^{ij}v_i)_j dt\) + \( -2 \sigma \ell_t v_t + 2\sum_{i,j=1}^n b^{ij}\ell_i v_j \)dt \\
      & \ds + \(\sigma (\ell_t ^2 - \ell_{tt} )  - \sum_{i,j=1}^n ( b^{ij}\ell_i \ell_j  - b^{ij}_j\ell_i  - b^{ij}\ell_{ij}) \) v dt\].
      \end{array}
\end{equation}
Therefore, by \eqref{4.15-eq1-eq2} and the
definition of $A$ in \eqref{AB1}, we
get
\begin{equation}\label{4.15-eq3}
\begin{array}{ll}
\ds\q\th\(-2\sigma\ell_tv_t+2\sum_{i,j=1}^n b^{ij}\ell_iv_j+\Psi v\)\( \sigma du_t-\sum_{i,j=1}^n(b^{ij}u_i)_jdt\)\\
\ns
 \ds= \(-2\sigma^2\ell_tv_t+2\sigma\sum_{i,j=1}^nb^{ij}\ell_iv_j+\sigma\Psi v\) dv_t\\
\ns
 \ds\q
 +\(-2\sigma\ell_tv_t+2\sum_{i,j=1}^n b^{ij}\ell_iv_j+\Psi v\)\(-\sum_{i,j=1}^n (b^{ij}v_i)_j+Av\)dt\\
\ns
 \ds\q+\(-2\sigma\ell_tv_t+2\sum_{i,j=1}^n b^{ij}\ell_iv_j+\Psi v\)^2dt.
\end{array}
\end{equation}

Let us continue to analyze the first two terms in
the right-hand side of \eqref{4.15-eq3}.

For the first term in the right-hand side of \eqref{4.15-eq3},  we find that
\begin{equation*}
\begin{cases}
  \ds\qquad -2\sigma^2 \ell_t v_t dv_t \3n&= \ds d(-\sigma^2 \ell_t v_t ^2) + \sigma^2 \ell_t (dv_t)^2 + (\sigma^2\ell_t )_t v_t^2 dt,\\
  \ds 2\sigma\sum_{i,j=1}^n b^{ij}\ell_i v_j d v_t \3n& =\ds  d\(2\sigma v_t \sum_{i,j=1}^n b^{ij}\ell_i v_j\) - 2\sum_{i,j=1}^n (\sigma b^{ij}\ell_i)_t v_j v_t dt- 2\sigma\sum_{i,j=1}^n b^{ij}\ell_i v_{jt}v_t dt,\\
  \ns\ds \qquad\qquad\sigma\Psi v dv_t \3n&\ds  =  d(\Psi \sigma v v_t) - (\sigma\Psi)_t v v_t dt - \sigma\Psi v_t^2 dt.
  \end{cases}
\end{equation*}
Therefore, we get that
\begin{equation}
\label{4.15-eq4}
  \begin{array}
    {ll}
    \quad &\ds \;\; \(-2\sigma\ell_tv_t+2\sum_{i,j=1}^nb^{ij}\ell_iv_j+\Psi v\)\sigma dv_t\\
 \ns & \ds =d\(-\sigma^2\ell_t v_t^2 +2\sigma v_t \sum_{i,j=1}^n b^{ij}\ell_i v_j + \sigma \Psi  v v_t - \frac 12 (\sigma\Psi)_t v^2\)\\
 \ns & \ds - \[\sum_{i,j=1}^n (\sigma b^{ij}\ell_i v_t^2)_j -\((\sigma\ell_t)_t + \sum_{i,j=1}^n (\sigma b^{ij}\ell_i)_j -\sigma \Psi\)v_t^2 \\
 & \ds + \; 2\sum_{i,j=1}^n (\sigma b^{ij}\ell_i)_t v_j v_t - \frac 12 (\sigma\Psi)_{tt} v^2\]dt + \sigma^2 \ell_t (d v_t)^2.
  \end{array}
\end{equation}

In a similar manner, for the second term in the right-hand side of \eqref{4.15-eq3},  we find that
\begin{equation}
\label{2nd in the right1}
\begin{array}{ll}\ds
 \quad \ds-2\sigma\ell _tv_t\[-\sum_{i,j=1}^n
(b^{ij}v_i)_j+Av\]\\
 \ns  = \ds 2\[\sum_{i,j = 1}^n
(\sigma b^{ij}\ell_tv_iv_t)_j-\sum_{i,j = 1}^n
 b^{ij}(\sigma\ell_{t})_j v_iv_t\]+\sum_{i,j = 1}^n
 (\sigma b^{ij}\ell_t)_tv_iv_j\\
\ns
  \ds\quad -\(\sigma\sum_{i,j = 1}^n
 b^{ij}\ell_tv_iv_j+\sigma A\ell_tv^2\)_t+(\sigma
 A\ell_t)_tv^2,
\end{array}
\end{equation}

\begin{equation}
\label{2nd in the right2}
\begin{array}{ll}\ds
 \quad   2\sum_{i,j=1}^n b^{ij}\ell_iv_j\[-\sum_{i,j = 1}^n  (b^{ij}v_i)_j+Av\]\\
\ns  =  \ds - \sum_{i,j= 1}^n \[\sum_{i',j'
=1}^n \(2b^{ij} b^{i'j'}\ell_{i'}v_iv_{j'}
-b^{ij}b^{i'j'}\ell_iv_{i'}v_{j'}\)-Ab^{ij}\ell_i v^2\]_j\\
\ns  \ds \q  +\sum_{i,j,i',j' = 1}^n
\left[2b^{ij'}(b^{i'j}\ell_{i'})_{j'} -
(b^{ij}b^{i'j'}\ell_{i'})_{j'}\right]v_iv_j-
\sum_{i,j = 1}^n (Ab^{ij}\ell_i)_j v^2,
\end{array}
\end{equation}
and
\begin{equation}
\label{2nd in the right3}
\begin{array}{ll}\ds
 \quad  \Psi v\[-\sum_{i,j=1}^n
(b^{ij}v_i)_j+Av\] &\ds = \ds  -\sum_{i,j = 1}^n
\(\Psi b^{ij}vv_i-\frac 12 {\Psi_i}b^{ij}
v^2\)_j+\Psi \sum_{i,j = 1}^n  b^{ij}v_iv_j \\
\ns &\quad \ds +\[-\frac 12 \sum_{i,j = 1}^n
(b^{ij}\Psi_i)_j+A\Psi\] v^2.
\end{array}
\end{equation}

Finally, from \eqref{4.15-eq3} to \eqref{2nd in
the right3}, we arrive at the desired equality
\eqref{hyperbolic2}.
\endpf



\section{Choice of the weight
function}\label{sec-weight}


In this section, we explain the choice of the
weight function which will be used to establish
our global Carleman estimate. Although such kind
of functions are already used in \cite{Amirov}, we give full details for the sake of
completion and the convenience of readers.

The weight function is given as follows:
\begin{equation}\label{weight1}
\varphi(x,t) = h x_1 + \frac{1}{2}\sum_{j=2}^n
x_j^2 + \frac{1}{2}\(t-\frac T2\)^2 +
\frac{1}{2}\tau,
\end{equation}
where $h$ and $\tau$ are suitable parameters, whose precise meanings will be explained in the sequel.

Without loss of generality, we assume that
$0=(0,\cds,0)\in S \setminus \partial S$ and
$\nu(0)=(1,\cds,0)$. For some $r >0$, for that $S$ is $C^2$, we can
parameterize $S$ in the neighborhood of the origin by
\begin{equation}\label{car eq1}
x_1=\g(x_2,\cds,x_n),\;|x_2|^2+\cds+|x_n|^2
<r.
\end{equation}
For natational brevity, denote $$a(x,t) = \frac{\partial \sigma}{\partial \nu}.$$
Hereafter, we set
\begin{equation}\label{car eq3}
\begin{cases}
\ds B_r\(0,\frac T2\)=\left\{(x,t):\,(x,t)\in\dbR^{n+1},\,|x|^2
+ \(t-\frac T2\)^2 < r^2\right\},\\[10pt]
 \ds B_r(0)=\{
x:\,x\in\dbR^n,\,|x|<r \}.
\end{cases}
\end{equation}
By \eqref{th1-eq1}, we have that
\begin{equation}\label{car eq2}
\left\{
\begin{array}{ll}\ds
-\a_0=a\(0,\frac T2 \)<0,\\[8pt]
\ns\ds
K <\frac{\a_0}{4(|\sigma|_{L^\infty(B_r(0,T/2))}+1)},\\[10pt]
\ns\ds - K \sum_{j=2}^n x_j^2
<\g(x_2,\cds,x_n),\;\mbox{ if
}\sum_{j=2}^n x_j^2 < r.
\end{array}
\right.
\end{equation}

Let
\begin{equation}\label{car eq4}
M_1=\max\left\{ |\sigma|_{C^1(B_r(0,0))},1 \right\}.
\end{equation}
Denote
$$
\cD^{-}_r=\{x:\,x\in
B_r(0),\,x_1<\gamma(x_2,\cds,x_n)\} ,\quad \cD^+_r
= B_r(0)\setminus \overline{D^-_r}.
$$
For any $\alpha \in (0, \alpha_0)$, in accordance with the continuity of $a(x,t)$ and the first
inequality in \eqref{car eq2}, it is clear that there exists
a $\delta_0>0$ small enough such that
$0<\d_0<\min\{1,r^2\}$, which would be specified later,  and
\begin{equation}\label{car eq5}
a(x,t)<-\alpha \;\mbox{ if } \;\;|x|^2 + \(t-\frac T2\)^2
\leq \d_0.
\end{equation}

Letting $ M_0=|\sigma|_{L^\infty(B_r(0,T/2))} $,
by the second inequality in \eqref{car eq2}, we
can always choose $K>0$ so large that
\begin{equation}\label{car eq6}
K<\frac{1}{2h}<\frac{\a}{4(M_0+1)}.
\end{equation}
Following immediately from \eqref{car eq6}, we have that
\begin{equation}
  \label{car eq6-cor}
  1-2hK > 0, \quad h\alpha - 2(M_0 + 1)> 0.
\end{equation}
For $K$ and $h$ such chosen, we will further
take  $\tau\in (0,1)$ so small that
\begin{equation}\label{car eq7}
\Big| \max\Big\{
\frac{K}{1-2hK},\frac{1}{2h} \Big\}
\Big|^2\tau^2 + \frac{2\tau}{1-2hK} \leq \delta_0.
\end{equation}
For convenience of notations, by denoting  $\mu_0 (\tau)$ the term in the left hand side of \eqref{car eq7} and letting $\cA_0 = \min\{\sigma, 1\}$,
we further assume that
\begin{equation}\label{car eq8-for M1 and N}
\left\{\begin{array}
{ll}
 h^2 \cA_0> 2hM_1 \sqrt{\mu_0(\tau)} + 2M_1\mu_0(\tau),\\[8pt]
 \ns\ds  \alpha h > 2(M_1^2 +
M_1)\sqrt{\mu_0(\tau)} - (M_0^2 + nM_0) -(n-1).
\end{array}
\right.
\end{equation}

For any positive number $\mu$ with $2\mu > \tau$, let
\begin{equation}\label{weight2}
Q_\mu = \left\{
(x,t)\in\dbR^{n+1}\| x_1>\gamma (x_2,x_3,\cdots,x_n), \sum_{j=2}^N x_j^2<\d_0,
\varphi(x,t)<\mu \right\}.
\end{equation}
The set $Q_{\tau}$ defined in this style is not empty. It is only to prove that the defining conditon $\varphi(x,t)<\mu$ is compatible with the first defining condition, i.e., $x_1>\gamma(x_2,x_3,\cdots,x_n)$.  By assumption, we know that $\gamma(x_2,x_3,\cdots, x_n) > -K \sum_{j=2}^n x_j^2$, then
together with the first inequality in \eqref{car eq6-cor}, we have that
\begin{eqnarray*}
  \varphi (x,t) \3n& \geq &\3n -hK \sum_{j=2}^n x_j^2 + \frac 12 \sum_{j=2}^n x_j^2 +\frac 12 \left(t-\frac T2\right)^2 +\frac 12 \tau\\
  & = &\3n \(\frac 12 -K h\)\sum_{j=2}^n x_j^2 +\frac 12 \(t-\frac T2\)^2 +\frac 12 \tau\\
  & > &\3n \frac\tau 2.
\end{eqnarray*}
Noting that  $(x,t)\in Q_{\mu}$ implies $\varphi(x,t)< \mu$, together with $\ds 2\mu > \tau$, we see by definition that
$Q_\mu \neq \emptyset$ as desired.

In what follows, we will show that how to determine the number $\delta_0$ appearing in \eqref{car eq7}.
Let $(x,t)\in \overline{Q}_{\tau}$. From the definition of $Q_{\tau}$ and noting that $\gamma(x_2,x_3,\cdots, x_n) > -K \sum_{j=2}^n x_j^2$, we find that
\begin{equation}\label{estimate x1-1}
x_1 \leq -\frac{\tau}{ 2h } \sum_{j=2}^n x_j^2 -
\frac 1{2h} \(t-\frac T2\)^2 + \frac{\tau}
{2h}\leq \frac\tau {2h}.
\end{equation}
On the other hand, by $\ds -K \sum_{j=2}^n x_j^2 \leq x_1$, we find that
\begin{equation*}
-K h \sum_{j=2}^n x_j^2 + \frac 12
\sum_{j=2}^n x_j^2 +\frac 12 \(t-\frac T2\)^2 +
\frac 12 \tau\leq \tau.
\end{equation*}
Thus
\begin{equation*}
  \sum_{j=2}^n x_j^2 < \frac\tau {1- 2K h}.
\end{equation*}
We then get that
\begin{equation}
\label{estimate x1-2}
-x_1\leq K \sum_{j=2}^n x_j^2< \frac{K\tau}{1 - 2K h}.
\end{equation}
Combining \eqref{estimate x1-1} and \eqref{estimate x1-2}, we arrive at
\begin{equation}
  \label{estimate x1-3}
  |x_1| \leq \max\left\{\frac{K }{1 -2hK}, \; \frac 1 {2h}\right\}\tau.
\end{equation}

Thus, by the restriction imposed on $\varphi (x,t)$ in the definition of $Q_{\tau}$ and \eqref{estimate x1-2}, we find that
\begin{equation}
\begin{array}{ll}
\ds \tau >   \ds  \varphi (x,t) = hx _1 + \frac 12
\sum_{j=2}^{n} x_j^2 + \frac 1 2 \(t-\frac T2\)^2 + \frac \tau 2\\
\ns \quad > \ds -\frac{K h \tau}{1- 2K
h} + \frac 12 \sum_{j=2}^n x_j^2 + \frac 12
\(t-\frac T2\)^2  +\frac \tau 2.
\end{array}
\end{equation}
This gives that
\begin{equation}
  \label{estimate t}
  \(t-\frac T2\)^2 <\frac{2K h \tau}{1 - 2K h}+ \tau = \frac\tau{1 - 2K h}.
\end{equation}
Correspondingly, we have that
\begin{equation*}
  |x|^2 + \(t-\frac T2\)^2 = x_1^2 + \sum_{j=2}^n x_j^2 + \(t-\frac T2\)^2 \leq  \|\max\left\{\frac K {1- 2K h}, \frac 1{2h}\right\}\|^2 \tau ^2 + \frac{2\tau} {1 -2K h}.
\end{equation*}

Returning back to \eqref{car eq5}, by \eqref{estimate x1-2}, \eqref{estimate x1-3} and \eqref{estimate t}, we choose the $\delta_0$ in the following style:
\begin{equation}
\label{choose delta0}
  \delta_0> \mu_0(\tau) =   \|\max\left\{\frac K {1- 2K h}, \frac 1{2h}\right\}\|^2 \tau ^2 + \frac{2\tau} {1 -2K h}.
\end{equation}
\ms


\section{A global Carleman
estimate}\label{sec-car}


This section is devoted to establishing a global
Carleman estimate for the stochastic hyperbolic
operator presented in Section 1, based on the
point-wise Carleman estimate given in Section \ref{sec-point}. It will be shown that it is
the key to the proof of the main result.

We have the following global Carleman estimate.
\begin{theorem}
\label{Car theorem} Let $u$ be an
$H^2_{loc}(\mathbb{R}^n)$-valued $\dbF$-adapted
process such that $u_t$ is an
$L^2(\mathbb{R}^n)$-valued semimartingale. If
$u$ is supported in $Q_{\tau}$, then there exist
a constant $C$ depending on $b_i, i=1,2,3$ and a
$s_0 > 0$ depending on $\sigma, \tau$ such that
for all $s \geq s_0$ it holds that
\begin{equation}\label{4.15-eq20}
\begin{array}{ll}\ds
\quad \mE\int_{Q_\tau}\theta \big( -2\sigma\ell_t v_t +
2\nabla\ell\cdot\nabla v \big)
\big( \sigma du_t - \Delta u dt \big)dx  \\
\ns\ds \geq C \mE\int_{Q_\tau} \[s\lambda^2
\varphi^{-\lambda-2}(|\nabla v|^2+v_t^2)+
s^3\lambda^4 \varphi^{-3\lambda-4}v^2 \]dxdt \\
\ns\ds \quad + \mE\int_{Q_\tau}\big( -2\sigma \ell_tv_t +
2\nabla\ell\cdot \nabla v \big)^2 dxdt + C
\mE\int_{Q_\tau}\sigma^2\theta^2\ell_t(du_t)^2dx.
\end{array}
\end{equation}
\end{theorem}

{\em Proof}. We apply the result of Lemma
\ref{hyperbolic1} to show our key Carleman
estimate. Let $(b^{ij})_{1\leq i,j\leq n}=I_n$ stand for
the unit matrix of $n$th order and let $\Psi = 0$ in \eqref{hyperbolic2}.
Then we find that
\begin{equation}\label{4.15-eq8}
\begin{array}{ll}
& \ds
\theta \left( -2\sigma\ell_t v_t +
2\n\ell\cdot\n v  \right)
\left[ \sigma du_t - \Delta u dt \right] \\[8pt]
\ns & \ds +\;\nabla\cdot\Big[ 2(\n
v\cd\n\ell)\n v - |\n v|^2\n \ell - 2
\ell_t v_t \n v
+ \sigma  v_t^2 \n \ell   -  A\nabla\ell\, v^2 \Big] dt \\[8pt]
\ns & \ds +\; d\Big[ \sigma \ell_t |\n v|^2 -
2\sigma \nabla\ell \cdot\nabla vv_t + \sigma^2\ell_t v_t^2 + \sigma A\ell_t  v^2 \Big] \\[8pt]
\ns\ds =  & \!\!\!\!\ds \Big\{ \[(\sigma^2\ell_{t})_t +
\nabla\cdot (\sigma \nabla \ell)  \]v_t^2\! -\! 2\[(\sigma
\n\ell)_t + \nabla (\sigma\ell_{t})\]\!\cdot\!\nabla v\; \! v_t +
\Big[ (\sigma \ell_t)_t + \Delta \ell
 \Big]|\nabla v|^2 \\[8pt]
\ns & \ds + Bv^2 + \Big(
-2\sigma\ell_tv_t + 2\n\ell \cdot\n v
\Big)^2\Big\} dt +
\sigma^2\theta^2\ell_t(du_t)^2,
\end{array}
\end{equation}
where $(du_t)^2$ represents the quadratic
variation process of $u_t$. It is easy to show that $(du_t)^2 = b_4^2 u^2 dt$ and $A$,
$B$ are stated  respectively as follows:
\begin{equation}\label{AB2}
\left\{
\begin{array}{ll} \ds
A\=\sigma(\ell_t^2-\ell_{tt})- (|\n\ell|^2
-\D \ell),\\
\ns \ds B\= (\sigma A\ell_t)_t- \nabla\cdot (A
\n\ell).
\end{array}
\right.
\end{equation}

Now let $\ell = s\varphi^{-\lambda}$ with $\varphi$  the weight function given by \eqref{weight1}. Then, some simple computations show that
\begin{equation}\label{4.15-eq9}
\left\{
\begin{array}{ll}\ds
\ell_t \3n& = \ds -s\lambda \varphi_t \varphi^{-n-1} = -s\lambda \(t-\dfrac T2\) \varphi^{-\lambda-1},\\[8pt]
\ns\ds \ell_{tt}\3n&\ds  = s\lambda (\lambda+1) \(t-\dfrac T2\)^2 \varphi^{-\lambda-2} - s\lambda\varphi^{-\lambda -1},\\[8pt]
\ns\ds \nabla \ell \3n&\ds = -s\lambda \varphi^{-\lambda -1}\nabla \varphi, \\[8pt]
\ns \Delta \ell \3n&\ds  = s\lambda (\lambda + 1)\varphi^{-\lambda -2}|\nabla \varphi|^2 - s\lambda \varphi^{-\lambda -1}\Delta \varphi,\\[8pt]
\nabla \ell_t  \3n& =\ds s\lambda (\lambda +
1)\varphi^{-\lambda -2} \(t-\dfrac T2\) \nabla
\varphi.
\end{array}
\right.
\end{equation}

We begin to analyze the terms in the right hand side of
\eqref{4.15-eq8} term by term. The first
one reads
\begin{equation}\label{4.15-eq10}
\begin{array}{ll}&\ds
\big[(\sigma^2\ell_{t})_t + \nabla\cdot(\sigma \n\ell)
\big]v_t^2 \\[8pt]
\ns\ds = &\3n \ds  \big[ 2\sigma\sigma_t \ell_t
+ \sigma^2 \ell_{tt} + \nabla
\sigma\cd\nabla\ell + \sigma\D\ell \big]v_t^2
\\[8pt]
\ns \ds = &\3n \ds\[2\sigma\sigma_t \ell_t + \sigma^2 \ell_{tt} -s\lambda (\nabla \sigma\cdot \nabla \varphi + \sigma\Delta\varphi)\varphi^{-\lambda -1} + s\lambda (\lambda +1)\sigma|\nabla \varphi|^2 \varphi^{-\lambda -2}\] v^2_t\\[8pt]
\ns =  &\3n\ds -s\lambda \varphi^{-\lambda -1}\[2 \sigma \sigma_t \;\(t-\dfrac T2\) + \sigma^2 +  (\nabla \sigma\cdot \nabla \varphi + \sigma\Delta \varphi) \] v_t ^2 \\[8pt]
\ns &\3n + s\lambda (\lambda +1)\varphi^{-\lambda -2}\[\sigma^2 \(t-\dfrac T2\)^2 + \sigma|\nabla \varphi|^2\] v_t ^2\\
\ns  \geq &\3n \ds - s\lambda \varphi^{-\lambda -1}\[h a + 2(M_1^2 +M_1)\sqrt{\mu_0(\tau)} + \(M_0^2 + (n-1)M_0\)\] v_t^2\\[8pt]
\ns  &\ds\3n + s\lambda (\lambda +1)h^2 \sigma  \varphi^{-\lambda -2} v_t ^2\\[8pt]
\ns  \geq &\3n \ds s\lambda \varphi^{-\lambda -1}\[h \alpha - 2(M_1^2 +M_1)\sqrt{\mu_0(\tau)} - \(M_0^2 + (n-1)M_0\)\] v_t^2\\[8pt]
\ns  &\3n\ds + h^2 \sigma  s\lambda (\lambda
+1)\varphi^{-\lambda -2} v_t ^2.
\end{array}
\end{equation}
Likewise, the second term in the right hand
side of \eqref{4.15-eq8} reads
\begin{equation}\label{4.15-eq11}
\begin{array}{ll}\ds
&\ds  - 2\[(\sigma \nabla \ell)_t + \nabla (\sigma\ell_{t})\]\cdot \nabla
v\,v_t \\[8pt]
\ns\ds = &\3n \ds - 2 \big[ \sigma_t\n\ell +
\sigma\nabla\ell_t + \ell_t\nabla \sigma +
\sigma\nabla \ell_t
\big]\cdot\nabla v \,v_t\\[8pt]
\ns = &\3n \ds \[2s\lambda \varphi^{-\lambda -1} (\sigma_t\nabla \varphi  + \(t-\dfrac T2\) \nabla \sigma) - 2s\lambda (\lambda + 1) \varphi^{-\lambda -2} \sigma t \nabla \varphi\]\cdot\nabla v\, v_t\\[8pt]
= &\3n \ds 2 s\lambda \varphi^{-\lambda -2}\[(\sigma_t \nabla \varphi +  \(t-\dfrac T2\)\nabla \sigma)\varphi - (\lambda + 1)\sigma \(t-\dfrac T2\)\nabla \varphi\]\cdot \nabla v \,\! v_t\\[8pt]
\ns \geq &\ds\3n -s \lambda \varphi^{-\lambda -2} \(M_1 h + 2 M_1\sqrt{\mu_0(\tau)}\)\tau \(|\nabla v|^2 + v_t ^2\)\\
\ns &\ds\3n - s\lambda (\lambda +
1)\varphi^{-\lambda -2}\(h M_1 \sqrt{\mu_0(\tau)} +
M_1\mu_0(\tau)\)\(|\nabla v|^2 + v_t ^2\).
\end{array}
\end{equation}
Thus, there exists a $\lambda_0 > 0$ such that for $\lambda > \lambda_0$ it holds that
\begin{equation}\label{4.15-eq11-1}
\begin{array}
{ll}
\quad - 2\[(\sigma\nabla \ell)_t + \nabla (\sigma\ell_{t})\]\cdot \nabla
v\,v_t\\
\ns  \geq \ds - 2s\lambda (\lambda + 1)\varphi^{-\lambda -2}\(h M_1 \sqrt{\mu_0(\tau)} + M_1\mu_0(\tau)\)\(|\nabla v|^2 + v_t ^2\).
\end{array}
\end{equation}

Treating the third term in the right hand side of \eqref{4.15-eq8} in the same fashion, we obtain
\begin{equation}\label{4.15-eq12}
\begin{array}{ll}\ds
& \ds\3n \big[(\sigma\ell_{t})_t + \Delta \ell
\big]|\nabla v|^2 \\[8pt]
\ns\ds = &\3n \ds \big[ \sigma_t \ell_t + \sigma
\ell_{tt} + \Delta\ell \big] |\nabla v|^2
\\[8pt]
\ns = & \3n\ds -s\lambda \varphi^{-n-1}\[\sigma_t \, \(t-\dfrac T2\) + \sigma + \Delta \varphi\]|\nabla v|^2 + \, s\lambda (\lambda +1) \varphi^{-\lambda -2}\[\sigma \(t -\dfrac T 2\)^2 +|\nabla \varphi|^2\]|\nabla v|^2\\
\ns \geq &\3n \ds -s\lambda \varphi^{-\lambda -
1}\(M_1 \sqrt{\mu_0(\tau)} + M_0 +
(n-1)\)|\nabla v|^2 + h^2 s\lambda (\lambda +
1)\varphi^{-\lambda -2}|\nabla v|^2.
\end{array}
\end{equation}

Following \eqref{4.15-eq10},  \eqref{4.15-eq11-1}, \eqref{4.15-eq12} and noticing \eqref{car eq8-for M1 and N}, we
find that
\begin{equation}\label{4.15-eq15}
\begin{array}{ll}\ds
& \ds \big[(\sigma^2\ell_{t})_t  + \nabla\cdot (\sigma \nabla\ell)
\big]v_t^2 - 2[(\sigma \nabla \ell)_t + \nabla (\sigma\ell_{t})]\cdot \nabla
v\, v_t  + \big[(\sigma\ell_{t})_t + \Delta\ell
\big]|\nabla v|^2 \\[8pt]
\ns &\ds\geq  C s\lambda^2
\varphi^{-\lambda-2}(|\nabla v|^2+v_t^2)
\end{array}
\end{equation}
for all $\lambda > \lambda_0$.

Next, note that in our case $A = \sigma(\ell_t ^2 - \ell_{tt}) - (|\nabla \ell |^2 - \Delta \ell)$. Then it is easy to show that
\begin{equation}\label{4.15-eq13}
\begin{array}{ll}\ds
A =& \3n \ds s^2 \lambda^2 \varphi^{-2\lambda -2}\[\sigma \(t -\dfrac T 2\)^2 - |\nabla \varphi|^2\] + s\lambda (\lambda +1)\varphi^{-\lambda -2} \[|\nabla \varphi|^2  - \sigma \(t -\dfrac T 2\)^2\]\\
\ns & \ds  + s\lambda \varphi^{-\lambda -1}[\sigma-
(n-1)].
\end{array}
\end{equation}
Thus, under some simple but a little more bothersome calculations, it holds that
\begin{equation}\label{4.15-eq14}
\begin{array}{ll}\ds
B \3n&\ds= (\sigma A\ell_t)_t- \nabla\cdot(A
\n\ell)\\[8pt]
\ns&\ds  = \sigma_t A\ell_t + \sigma
A_t\ell_t+\sigma A\ell_{tt}-\nabla A\cd\nabla\ell - A\D\ell \\[8pt]
\ns&\ds = 3s^3\lambda^2(\lambda+1)^2
\(t -\dfrac T 2\)^2\[\(t -\dfrac T 2\)^2-|\nabla \varphi|^2\]\varphi^{-3\lambda-4}\\[8pt]
\ns&\ds \q + 3s^3\lambda^2(\lambda+1)^2
|\nabla \varphi|^2\[|\nabla \varphi|^2-\(t -\dfrac T 2\)^2\]\varphi^{-3\lambda-4}\\[8pt]
\ns&\ds \q + O(s^3\lambda^3\varphi^{-3\lambda-3}) +
O(s^2\lambda^4\varphi^{-3\lambda-4}).
\end{array}
\end{equation}
It is easy to see that there exist an $\lambda_1>0$ and $s_0>0$ such
that for all $\lambda\geq \lambda_1$, $s\geq s_0$,
\begin{equation}\label{4.15-eq17}
\begin{array}{ll}\ds
Bv^2 \geq Cs^3\lambda^4
\varphi^{-3\lambda-4}v^2.
\end{array}
\end{equation}

Next, integrating \eqref{4.15-eq8} over $Q_\tau$ and
taking mathematical expectation, we obtain that
\begin{equation}\label{4.15-eq18}
\begin{array}{ll}\ds
\q \mE\int_{Q_\tau}\theta \big( -2\sigma\ell_t v_t +
2\nabla\ell\cdot\nabla v \big)
\big( \sigma du_t - \D u dt \big)dx  \\
\ns\ds \geq C \mE\int_{Q_\tau} \[s\lambda^2
\varphi^{-\lambda-2}(|\nabla v|^2+v_t^2)+
s^3\lambda^4 \varphi^{-3\lambda-4}v^2 \]dxdt \\
\ns\ds \q + \mE\int_{Q_\tau}\big( -2\sigma\ell_tv_t +
2\nabla\ell \cdot\nabla v \big)^2 dxdt + C
\mE\int_{Q_\tau}\sigma^2\theta^2\ell_t(du_t)^2dx.
\end{array}
\end{equation}
Thus we complete the proof.\endpf

\section{Proofs of the Main
Result}\label{sec-main}

This section is dedicated to the proof of the unique continuation property presented in Section 1.

{\em Proof}. Without loss of generality, we assume that $$x_0 = (0,0,\cdots, 0),\quad \nu (x_0) = (1,0,\cdots, 0)$$ and $S$ is parameterized as in Section \ref{sec-weight} near $0$. Also, $K, \delta_0, h, \tau$ are all given as in Section \ref{sec-weight}. By the definition of $\varphi(x,t)$ and $Q_{\mu}$, for any $\mu\in (0, \tau]$, the boundary $\Gamma_\mu$ of $Q_{\mu}$ is composed of the following three parts:
\begin{equation}
  \label{boundary of the domain}
 \hspace{-0.5cm}\begin{cases}
    \ds \Gamma_\mu ^1 = &\hspace{-0.3cm} \ds \bigg\{(x,t) \in \dbR^{n+1} \| x_1 = \gamma(x_2, x_3, \cdots, x_n), \sum_{j=2}^n x_j^2 < \delta_0, \varphi(x,t)<\mu \bigg\},\\[12pt]
     \ds \Gamma_\mu ^2 = & \hspace{-0.3cm}\ds \bigg\{(x,t) \in \dbR^{n+1} \| x_1 > \gamma(x_2, x_3, \cdots, x_n), \sum_{j=2}^n x_j^2 < \delta_0, \varphi(x,t)=\mu \bigg\},\\[12pt]
      \ds \Gamma_\mu^3 = &\hspace{-0.3cm} \ds \bigg\{(x,t) \in \dbR^{n+1} \| x_1 > \gamma(x_2, x_3, \cdots, x_n), \sum_{j=2}^n x_j^2 = \delta_0, \varphi(x,t)<\mu \bigg\},
  \end{cases}
\end{equation}
i.e., $\ds \Gamma_{\mu} =\Gamma_\mu^1\cup \Gamma_\mu^2\cup \Gamma_\mu^3$. Next, we show that in fact $\Gamma_\mu^3 = \emptyset$. Based on the conditions $\ds x_1 > \gamma (x_2, x_3, \cdots, x_n)$ and $\gamma(x_2,x_3,\cdots, x_n) >-K \sum_{j=2}^n x_j^2$ and the definition of $\varphi$, it follows that
\begin{equation}
\label{Gamma 3 empty}
  (1-2K h) \sum_{j=2}^n x_j^2 + \(t -\dfrac T 2\)^2 < 2\[h x_1 + \sum_{j=2}^n x_j^2 + \(t -\dfrac T 2\)^2\]= 2\varphi - \tau < 2\mu - \tau <\tau.
\end{equation}
Also, note that $\Gamma_\mu^3$ is subordinated to $\sum_{j=2}^n x_j^2 = \delta_0$. Thus, from \eqref{Gamma 3 empty}, it follows that $\delta_0 < \frac{\tau}{1 - 2K h}$, a contradiction to $\delta_0 > \frac{\tau}{1 -2K h}$ introduced in \eqref{choose delta0}. As a direct result, we conclude that $\Gamma_{\mu} = \Gamma_\mu^1\cup \Gamma_\mu^2$.

Moreover, it is clear that
$$\Gamma_\mu^1\cup \Gamma_\mu^2 \subset \overline{Q_{\tau}}.$$
Define $$t_0 = \sqrt{\frac{\tau}{1  - 2K h}}\;.$$
Then by \eqref{estimate t}, it follows
\begin{equation}
  \label{boundary containing}
  \begin{cases}
    \Gamma_\mu^1 \subset \{x\mid x _1 = \gamma (x_2, x_3, \cdots, x_n)\} \times \left\{t \mid \;|\,t-T/2\,| \leq t_0\right\},\\
    \ns \Gamma_\mu^2 \subset \{x\mid \varphi(x,t) = \mu\}, \quad \mu \in (0, \tau].
  \end{cases}
\end{equation}
It is clear that
\begin{equation*}
  \Gamma_\mu^j \subset \Gamma_\tau^j, \quad j =1,2.
\end{equation*}

To apply the result of Theorem \ref{Car theorem} to the present case, we adopt the truncation method. For convenience in the later statement, denote $Q_{\tau} = Q_1$. Fixing a arbitrarily small number $\widetilde \tau \in (0, \frac{\tau}{8})$, let
\begin{equation*}
  Q_{k+1} = \{(t,x)| \varphi(x,t)<\tau - k\widetilde \tau, k=1,2,3\}.
\end{equation*}
Hence, it is easy to show that that $Q_4\subset Q_3\subset Q_2\subset Q_1$.

Introduce a truncation function $\chi\in C_0^{\infty}(Q_2)$ in the following manner
\begin{equation*}
  \chi \in [0, 1]\quad \text{and}\quad \chi = 1 \quad \text{in} \quad Q_3.
\end{equation*}

Let $z$ be the solution of \eqref{system1}. Let $\Phi= \chi z$. Then a little bothersome calculation gives that
\begin{equation}
\label{truncated equation}
\begin{cases}
\sigma d\Phi - \Delta \Phi dt = (b_1\Phi_t + b_2\cdot\nabla \Phi + b_3 \Phi)dt  + f(x,t)dt +  b_4\Phi dW(t)& \text{in } Q_{\tau}, \\
\ns \ds \Phi = 0, \;\; \frac{\partial
\Phi}{\partial \nu} = 0 &\text{on }
\Gamma_{\tau}.
\end{cases}
\end{equation}
Here, we denote by
\begin{equation*}
  f(x,t) = \sigma \chi_{tt} z + 2a\chi_t z_t - 2\nabla \chi\cdot\nabla z - z\Delta \chi - b_1\chi_t z - b_2\cdot z\nabla \chi.
\end{equation*}
From the definition of $\chi$, $f$ is clearly supported in $Q_2\setminus \overline Q_3$.

Let $$F = b_1\Phi_t + b_2\cdot \nabla \Phi + b_3 \Phi + f.$$
In stead of $u$ by $\Phi$ in \eqref{4.15-eq20}, we have that
\begin{equation}\label{4.15-eq18}
\begin{array}{ll}\ds
\q \mE\int_{Q_\e}\theta \big( -2\sigma\ell_t v_t +
2\nabla\ell\cdot\nabla v \big)
\big( \sigma d\Phi_t - \Delta \Phi dt \big)dx  \\
\ns\ds \geq C \mE\int_{Q_\e} \[s\lambda^2
\varphi^{-\lambda-2}(|\nabla v|^2+v_t^2)+
s^3\lambda^4 \varphi^{-3\lambda-4}v^2 \]dxdt \\
\ns\ds \quad + \mE\int_{Q_\tau}\big( -2\sigma\ell_t v_t +
2\nabla\ell \cdot\nabla v \big)^2 dxdt + C
\mE\int_{Q_\tau}\sigma^2\theta^2\ell_t(d\Phi_t)^2dx.
\end{array}
\end{equation}
Due to the elementary property of It\^{o} integral, it is clear that
 \begin{equation*}
   \begin{array}
     {ll}
     & \quad\ds \mE\int_{Q_\tau}\theta \big( -2\sigma\ell_t v_t +
2\nabla\ell\cdot\nabla v \big)
\big( \sigma d\Phi_t - \Delta \Phi dt \big)dx\\
\ns & =\ds \mE\int_{Q_\tau}\theta \big( -2\sigma\ell_t v_t +
2\nabla\ell\cdot\nabla v \big)F dxdt \\
\ns & \ds\quad  + \; \mE\int_{Q_\tau}\theta \big( -2\sigma\ell_t v_t +
2\nabla\ell\cdot\nabla v \big)b_4 \Phi d W(t) dx\\
\ns & \leq \ds  \mE \int_{Q_{\tau}}\theta^2 F^2 dxdt + \mE\int_{Q_\e}\theta \big( -2\sigma\ell_t v_t +
2\nabla\ell\cdot\nabla v \big)^2 dxdt.
   \end{array}
 \end{equation*}
 Thus, we show that
 \begin{equation}
   \begin{array}
     {ll}
    &\ds \mE\int_{Q_\tau} \[s\lambda^2
\varphi^{-\lambda-2}(|\nabla v|^2+v_t^2)+
s^3\lambda^4 \varphi^{-3\lambda-4}v^2 \]dxdt \\
\ns & \ds +\;
\mE\int_{Q_\tau}\sigma^2\theta^2\ell_t(d\Phi_t)^2dx \leq C\mE \int_{Q_{\tau}}\theta^2 F^2 dxdt.
   \end{array}
 \end{equation}
 Let us now do some estimate for the right hand side of the above inequality.
 \begin{equation}
   \begin{array}
     {ll}
     & \quad \ds \mE \int_{Q_{\tau}}\theta^2 F^2 dxdt\\
     \ns & \ds\leq 2\mE\int_{Q_{\tau}}\(b_1\Phi_t + b_2\cdot\nabla \Phi + b_3 \Phi\)^2 dxdt + 2\mE\int_{Q_{\tau}}|f|^2 dxdt.
   \end{array}
 \end{equation}
 Note that $f$ is supported in $Q_2\setminus\overline Q_3$. Hence
 \begin{equation}
   \begin{array}
     {ll}
     & \ds \quad \mE\int_{Q_{\tau}}\theta^2|f|^2 dxdt\\
     \ns & \ds = \mE\int_{Q_{\tau}}\theta^2| \sigma \chi_{tt} z + 2\sigma\chi_t z_t - 2\nabla \chi\cdot\nabla z - z\Delta \chi - b_1\chi_t z - b_2\cdot z\nabla \chi|^2 dxdt\\
     \ns & \ds \leq C \mE\int_{Q_2\setminus\overline Q_3} \theta^2\(z_t^2 + |\nabla z|^2 + z^2\)dxdt.
   \end{array}
 \end{equation}
 Thus, we achieve that
 \begin{equation}
 \begin{array}
 {ll}
\q \ds \mE \int_{Q_{\tau}}\theta^2 F^2 dxdt\\
   \ns   \ds \leq C \mE \int_{Q_1}\theta^2 \(\Phi_t ^2 + |\nabla \Phi|^2 + \Phi^2\)dxdt+  \; C \mE\int_{Q_2\setminus \overline Q_3}
\theta^2\(z_t^2 + |\nabla z|^2 + z^2\)dxdt.
   \end{array}
 \end{equation}
  And then
  \begin{equation}
   \begin{array}
     {ll}
\q\ds \mE\int_{Q_\tau} \[s\lambda^2
\varphi^{-\lambda-2}(|\nabla v|^2+v_t^2)+ s^3\lambda^4
\varphi^{-3\lambda-4}v^2 \]dxdt  +
\mE\int_{Q_\e}\sigma^2\theta^2\ell_t(d\Phi_t)^2dx \\
   \ns   \ds \leq C \mE \int_{Q_1}\theta^2 \(\Phi_t ^2 + |\nabla \Phi|^2 + \Phi^2\)dxdt+  \; C \mE\int_{Q_2\setminus \overline Q_3}
\theta^2\(z_t^2 + |\nabla z|^2 + z^2\)dxdt.
   \end{array}
 \end{equation}
 Recall that
 \begin{equation*}
   (d\Phi_t)^2  = b_4^2 \Phi^2dt, \quad \ell_t = -s\lambda \(t-\dfrac T2\) \varphi^{-\lambda -1}.
 \end{equation*}
 Therefore
 \begin{equation}
   \hspace{-0.5cm}\begin{array}
     {ll}
    &\ds \mE\int_{Q_{\tau}} \[s\lambda^2
\varphi^{-\lambda-2}(|\nabla v|^2+v_t^2)+
s^3\lambda^4 \varphi^{-3\lambda-4}v^2 \]dxdt  \\
   \ns & \ds \leq C \mE \int_{Q_1}\theta^2 \(\Phi_t ^2 + |\nabla \Phi|^2 + \Phi^2\)dxdt +
\mE\int_{Q_\tau}\theta^2 \sigma^2 s\lambda t \varphi^{-\lambda - 1} b_4^2\Phi ^2dxdt\\
    \ns & \ds \quad +  \; C \mE\int_{Q_2\setminus \overline Q_3} \theta^2\(z_t^2 + |\nabla z|^2 + z^2\)dxdt.
   \end{array}
 \end{equation}

 Then for $s$ and $\lambda$ large enough, it follows that

\begin{equation}
  \begin{array}
     {ll}
    &\ds \mE\int_{Q_{\tau}} \[s\lambda^2
\varphi^{-\lambda-2}(|\nabla v|^2+v_t^2)+
s^3\lambda^4 \varphi^{-3\lambda-4}v^2 \]dxdt  \\
   \ns & \ds \leq C \mE \int_{Q_\tau}\theta^2 \(\Phi_t ^2 + |\nabla \Phi|^2 + \Phi^2\)dxdt  +  \; C \mE\int_{Q_2\setminus\overline Q_3} \theta^2\(z_t^2 + |\nabla z|^2 + z^2\)dxdt.
   \end{array}
 \end{equation}

 Now, we find that
 \begin{equation*}
   |\nabla v|^2 + v_t^2 \geq C \theta^2 \(s^2\lambda^2 \varphi^{-2\lambda  -2}\Phi^2 + |\nabla \Phi|^2 + \Phi_t^2\).
 \end{equation*}
 Thus for large $s$ and $\lambda$, it follows that
 \begin{equation}
   \begin{array}
     {ll}
    &\ds \mE\int_{Q_{\tau}}\theta^2 \[s\lambda^2
\varphi^{-\lambda-2}(|\nabla \Phi|^2+\Phi_t^2)+
s^3\lambda^4 \varphi^{-3\lambda-4}\Phi^2 \]dxdt  \\
   \ns & \ds \leq  C \mE\int_{Q_2\setminus\overline Q_3} \theta^2\(z_t^2 + |\nabla z|^2 + z^2\)dxdt.
   \end{array}
 \end{equation}
Recall that  and $\Phi = z$ in $Q_3\subset Q_{\tau}$. It is easy to show that
\begin{equation}
\label{reduce to subinterval-1}
  \begin{array}
     {ll}
    &\ds \mE\int_{Q_{3}}\theta^2 \[s\lambda^2
\varphi^{-\lambda-2}(|\nabla z|^2+z_t^2)+
s^3\lambda^4 \varphi^{-3\lambda-4}z^2 \]dxdt  \\
   \ns & \ds \leq  C \mE\int_{Q_2\setminus \overline Q_3} \theta^2\(z_t^2 + |\nabla z|^2 + z^2\)dxdt.
   \end{array}
 \end{equation}

 Note that in $Q_4$, $\ds \varphi(x,t) < \tau - 3\widetilde \tau$, then $\theta = e^{s\varphi^{-\lambda}} > e^{s(\tau - 3\widetilde \tau)^{-\lambda}}$. Moreover, in $Q_2\setminus \overline Q_3$, $\ds \tau - 2\widetilde \tau < \varphi(x,t) <\tau -\widetilde \tau$, then $e^{s(\tau - \widetilde \tau)^{-\lambda}} < \theta < e^{s(\tau - 2\widetilde \tau)^{-\lambda}}$. Therefore
 \begin{equation}
\label{reduce to subinterval-2}
  \begin{array}
     {ll}
    &\ds \mE\int_{Q_{4}}\[s\lambda^2
\varphi^{-\lambda-2}(|\nabla z|^2+z_t^2)+
s^3\lambda^4 \varphi^{-3\lambda-4}z^2 \]dxdt  \\
   \ns & \ds \leq  C e^{2[s(\tau-2\widetilde \tau)^{-\lambda} - s(\tau - 3\widetilde \tau)^{-\lambda}]}\mE\int_{Q_2\setminus \overline Q_3} \(z_t^2 + |\nabla z|^2 + z^2\)dxdt\\
   \ns & \ds \leq  C e^{2[s(\tau-2\widetilde \tau)^{-\lambda} - s(\tau - 3\widetilde \tau)^{-\lambda}]}\mE\int_{Q_\tau} \(z_t^2 + |\nabla z|^2 + z^2\)dxdt
   \end{array}
 \end{equation}
 For brevity, by letting $\overline \mu = 2[(\tau-2\widetilde \tau)^{-\lambda} - (\tau - 3\widetilde \tau)^{-\lambda}]$, we can have that
 \begin{equation}
\label{reduce to subinterval-2}
  \begin{array}
     {ll}
    &\ds \mE\int_{Q_{4}}(|\nabla z|^2+z_t^2+ z^2 )dxdt  \leq  C e^{s\overline \mu} \mE\int_{Q_\tau} \(z_t^2 + |\nabla z|^2 + z^2\)dxdt.
   \end{array}
 \end{equation}
 For that $\overline \mu < 0$, so if we let $s\to +\infty$, we find $z = 0$ in $Q_4$. Taking $Q_4$ the desired region, we complete the proof. \endpf




\end{document}